\documentclass[11pt, letter, oneside]{article}

\usepackage{epsfig}
\usepackage{fancyhdr}
\usepackage{amssymb}
\usepackage[mathscr]{eucal}
\usepackage[centertags]{amsmath}
\usepackage{newlfont,color}
\usepackage{latexsym,enumerate,epsfig,listings}
\usepackage{amsmath,amsthm,amsopn,amstext,amscd,amsfonts,amssymb}
\usepackage{graphicx}
\usepackage{hyperref}
\usepackage{float}
\usepackage[sectionbib]{chapterbib} 

\newtheorem{teo}{Theorem}
\newtheorem{lema}{Lemma}

\def \beq {\begin{eqnarray}}
\def \eeq {\end{eqnarray}}
\def \beqn {\begin{eqnarray*}}
\def \eeqn {\end{eqnarray*}}

\newcommand{\beeq}{\begin{equation}}
\newcommand{\eneq}{\end{equation}}
\newcommand{\bearno}{\begin{eqnarray*}}
\newcommand{\enarno}{\end{eqnarray*}}
\newcommand{\befi}{\begin{figure}}
\newcommand{\enfi}{\end{figure}}

\def\ep{\hfill $\Box$}

\def\bp{\noindent{\bf Proof.}\ }

\def\Q{{\cal Q}}

\newcommand{\pro}{\mathbb{P}}

\newcommand{\esp}{\mathbb{E}}

\newcommand{\na}{ \mathbb{N} }

\begin{document}

\author{Sergio I. L\'opez}
\date{}
\title{Convergence of tandem Brownian queues.}

\maketitle

\begin{abstract}
It is known that in a stationary Brownian queue with both arrival and service processes equal in law to Brownian motion, the departure process is a Brownian motion, that is, Burke's theorem in this context. In this short note
we prove convergence to this invariant measure: if we have an arbitrary continuous process satisfying some mild conditions as initial arrival process and pass it through an infinite tandem network of queues, the resulting process weakly converges to a Brownian motion. We assume independent and exponential initial workloads for all queues. 
\end{abstract}

\section{Introduction}\label{sec:intro1}

In 1956 Burke \cite{Burke} obtained one fundamental result for queueing theory. The first part of this result states that given an arrival process Poisson with rate $\lambda < 1$, and an independent 
service Poisson process with rate $1$, which together define a $M/M/1$ queue, the departure process is Poisson 
with parameter $\lambda$. The second part states a factorization property: the length of the queue at time $t$ is independent of future arrivals and past departures. Several extensions of this result have followed, see for example
\cite{Draief2003, Ferrari94, Martin10}.

The Brownian queue is a continuous valued model for a queue, which is indeed the heavy traffic limit of a $M \setminus M \setminus 1$ queue. We define it in the following. Denote by $\mathcal R: D[0, \infty) \rightarrow D[0, \infty)$ the operator in the space of $c\grave{a}dl\grave{a}g$ functions given by
\begin{equation}
\mathcal R (f) (t) :=  f(t) - \inf_{0 \leq s \leq t}  \{ f(s) \wedge 0 \},
\end{equation} 
called Skorokhod reflective mapping, see \cite{Harrison90b} pp. 14, for instance. Then, given two functions $f,g \in D[0, \infty)$ such that $f(0) \geq g(0)$, one can define the reflection of one function on another one by next mapping:
\begin{equation}\label{birefl}
L_f(g) \, (t) :=  g(t) + \inf _{0 \leq s \leq t} \Big\{ (f(s)-g(s)) \wedge 0 \Big\} = f- \mathcal R(f-g).
\end{equation}
Heuristically, the function $L_f(g)$ is such that has $g$ as driving function and $f$ as an upper barrier, where
non-ellastic collisions take place. This mapping provides the definition of a queue in the context of $c\grave{a}dl\grave{a}g$ functions by next interpretation: $f$ is the arrival process, $g$ the service process, and $f(0)-g(0)$
the initial workload of the queue. The queueing operator is then defined by 
\begin{equation}\label{queue}
\mathcal{Q} (f,g): = L_f(g),
\end{equation}
meaning that $\mathcal{Q} (f,g)$ is the departure process. Other important processes are the queue length process given by $f-\mathcal{Q} (f,g)$ and the free process defined as $f-g$. 


Take the particular case where $f$ is drived by a standard Brownian motion, $g$ drived by a Brownian motion with positive drift $c$ and choose the initial workload $f(0)-g(0)$ as the stationary distribution of 
this queue (namely an exponential random variable of parameter $c$, see \cite{Harrison90b}) pp. 15, all independently. Then, the above definition matches the one given by O'Connell and Yor \cite{OConnell01} of the stationary version of the 
Brownian queue (for positive times), as Norros and Salminen \cite{Salminen01} pointed out. 

For this model (and further generalizations of functionals of Brownian motion) an analogue result to Burke's theorem is presented by O'Connell and Yor in \cite{OConnell01}:
\begin{teo}\label{BROW:B}
Let $B^1_t$, $B^2_t$ be standard Brownian motions, $\mathcal{E}$ an exponential variable of parameter $c$, and $x$ a real number. Assume that all random elements are independent. Then
\begin{enumerate}
\item $ D_t = \mathcal{Q}( B^1_t + x, B^2_t + c\, t + x - \mathcal{E})$ has the law of a standard Brownian motion starting at $x -\mathcal{E}$,
\item $\{ D_s :  0 \leq s < t \}$ and $Q_t$ are independent.
\end{enumerate}
\end{teo}
The proof of this result goes backs to Harrison and Williams \cite{Harrison90}, in the context of multiclass stations, and it relies on weak convergence arguments, or alternatively, on path properties of the Brownian motion.

A tandem queue is a system of queues where there is an arrival process $A^1$, and a sequence $\{ S^n \}_{n \geq 1}$ of service processes, all independent. The system is defined recursively. The initial queue is fed from the arrival process $A^1$, and has departures determined by the service process $S^1$. For $n \geq 2$, the arrival process for the $n$-th queue is defined as the departure process of the $(n-1)$-th queue and the departures are determined by the service process $S^n$. 

When the initial arrival process is Poisson, Burke's theorem allows us to treat a tandem system of queues at any fixed time as if the queues acted independently. For example, take a two nodes system of tandem queues, with arrivals Poisson($\lambda$), $\lambda <1$, service processes Poisson($1$), all independent, and sample the initial length of each queue from its stationary measure. 
Because of the first part of Burke's theorem, the departure process of the first node is a Poisson($\lambda$) process. Moreover, due to the second part of Burke's theorem, the departure process of the first queue prior to time $t$ is independent of $Q^1_t$, the length of that queue at time $t$. Then the length of the second queue at time $t$, $Q^2_t$, is independent from $Q^1_t$, and it follows that the invariant measure of the system is a product measure. The factorization property from Burke's theorem has thus enabled the analysis of more complex systems. 


In the case when the initial arrival process is not Poisson, a natural question is whether it is possible to prove convergence to the stationary distribution in a tandem system where the number of queues goes to infinite. Assuming an existence result, Anantharam \cite{anan} proved the uniqueness of a stationary ergodic fixed point for the $\cdot$/M/K queue. Next, Mountford and Prabhakar \cite{MountfordPrabhakar} proved the attractiveness of the Poisson distribution in the class of ergodic stationary point processes on the line. To obtain this result, they used a coloring coupling technique based on an argument of Ekhaus and Gray (unpublished, cited by \cite{MountfordPrabhakar}).

In this note we present an analogue of Mountford-Prabhakar's theorem, for the Brownian queue:
\begin{teo}\label{BROW:MP}
Let $A^0$ be a process with continuous paths $A^0(\cdot,\omega):[0,\infty) \rightarrow \mathbb{R}$ that do not explode in finite time almost surely, and $A^0(0,\omega) \equiv 0$. Let $\{ W^ n \}_{n \in \na}$ be a family of standard Brownian motions, $\{ \mathcal{E}^n \}_{n \in \na}$ a family of exponential random variables with common parameter $c>0$, all independent. We define recursively the sequence of processes
$$  A^n = \Q \Big( A^{n-1}, W^n + c \, t - \sum_{i=1}^{n} \mathcal{E}^i \Big) , \quad  n \geq 1,$$
where $\Q$ is the queueing operator defined in Equation \ref{queue}. Then $A^n + \sum_{i=1}^{n} \mathcal{E}^i $ weakly converges to a Brownian motion.
\end{teo}

In words, the departure process of a tandem system of Brownian queues is weakly convergent to Brownian motion, for an initial arrival process belonging to a wide class of continuous-valued processes, and a particular set of initial conditions for the tandem queues: all having independent workloads, distributed as the stationary distribution of the Brownian queue. 
The coupling used in the $\cdot \setminus M \setminus 1$ case \cite{MountfordPrabhakar} is no longer suitable and we 
introduce an ad-hoc coupling technique that takes advantage of simple path properties of the Brownian motion. 
This procedure, however, strongly depends on the particular choice of initial workloads for the queues. 
We are currently working on a version of Theorem \ref{BROW:MP} where each queue is stationary, using a different
approach.  

For completeness, we present an elementary proof of Theorem \ref{BROW:B}, using the heavy traffic weak limit of the 
$M \setminus M \setminus 1$ queue, as was done in \cite{Harrison90b}, but avoiding the more complex context of 
multiclass stations.

\section{Proofs}\label{Sec:BA}

\subsection{Burke's Theorem analogue}

Before proving Theorem \ref{BROW:B} we state a corollary of Donsker's theorem:

\begin{lema}\label{CLPoi}
Let $\{ P^n \}_{n \in \na}$ be a sequence of Poisson processes with rate $r_n>0$. Suppose that $r_n \rightarrow r \in (0, \infty)$. Then
$$ \frac{ P^n(nt) - r_n \, n \, t }{ \sqrt{n} }  \Rightarrow \sqrt{r} \, B(t),$$
where $\{ B(t) \}_{ t \geq 0 }$ is a standard Brownian motion, and $\Rightarrow$ denotes the weak convergence of processes.
\end{lema}

\bp Let $P(t)$ be a Poisson process with intensity $1$. Let us denote equality in law by $=^d$. Since $ \{ P^n (t) : t \geq 0 \}  =^d \{ P ( r_n \, t ) : t \geq 0 \}$, we have
$$ \frac{P^n(nt) - n \, r_n t}{\sqrt{n}} =^d \frac{ P(n \,r_n \, t ) - n \, r_n \, t}{\sqrt{n}} =
\Big( \frac{ P( (n \,r_n) \, t ) - (n \, r_n) \, t}{\sqrt{n \, r_n}} \Big) \, \sqrt{r_n} .$$
The result follows by the functional CLT. \ep

\textbf{Proof of Theorem \ref{BROW:B}:} 

$1$) Define $\lambda_n := 1 - \frac{c}{\sqrt{n}}$. For each $n$ large enough such that $0 < \lambda_n <1$, we define $A^n (t)$ be a Poisson process with parameter $\lambda_n$, $S(t)$ a Poisson process with parameter $1$, 
and $G^n$ a geometric random variable with parameter $ \lambda_n$, all independent. Next, let 
$$ \tilde{A}^n(t) :=  \frac{ A^n(nt)-n  t \lambda_n }{ \sqrt{n} }. $$
Applying Lemma \ref{CLPoi}, $ \tilde{A}^n(t)$ weakly converges to a Brownian motion. We also define 
$$ \tilde{S}^n(t) := \frac{S(nt) - n t \lambda_n}{\sqrt{n}}. $$
We have that $\sqrt{n} (1- \lambda_n)=c$ and then
$$ \tilde{S}^n(t) = \frac{S(nt) - n t }{\sqrt{n}} + c \,t  ,   $$
and this weakly converges to $B(t) + c \, t$, where $B(t)$ is a Brownian motion, by Lemma \ref{CLPoi}. Finally, let $\tilde{G}^n := \frac{G^n}{\sqrt{n}}$, so that $\tilde{G}^n$ converges to an exponential random variable of parameter $c$.

Since $A^n, S^n$ and $G^n$ are independent, we conclude that $(\tilde{A}^n(t) + x,\tilde{S}^n(t) + x - \tilde{G}^n )$ weakly converges to $(B^1(t)+x,B^2(t) +c\, t + x - \mathcal{E})$, where $B^1$ and $B^2$ are standard Brownian motions
and $\mathcal{E}$ is an exponencial random variable of parameter $c$, all independent. The reflective mapping is continuous in the Skorokhod topology (\cite{Whitt} pp. 439) and then the queueing operator is 
continuous. By the continuous mapping theorem we have that $\mathcal Q (\tilde{A}^n(t)+x, \tilde{S}^n(t) + x - \tilde{G}^n )$ weakly converges to $\mathcal Q(B^1(t)+x,B^2(t) +c\, t + x -\mathcal{E})$. 

On the other hand, 
let $D^n(t) := \mathcal Q( A^n(t), S(t) - G^n)$. By Burke's theorem, $D^n$ has the law of a rate $ \lambda_n$ Poisson process. Define the process
$$ \tilde{D}^n(t) := \frac{ D^n(nt)- n \, t \, \lambda_n }{ \sqrt{n} }. $$
Again, By Lemma \ref{CLPoi}, $\tilde{D}^n$ weakly converges to a standard Brownian motion $B(t)$. Note that
\begin{eqnarray*}
 \tilde{D}^n(t) & = & \frac{ \mathcal Q( A^n(nt), S(nt) - G_n) - n \, t \, \lambda_n }{\sqrt{n}} \\
& = &  \mathcal Q \Big( \frac{ A^n(nt) - n \, t \, \lambda_n}{\sqrt{n}} , \frac{ S(nt) - n \, t \, \lambda_n - G_n}{\sqrt{n}} \Big) = \mathcal Q ( \tilde{A}^n(t),\tilde{S}^n(t) - \tilde{G}^n )  ,
\end{eqnarray*}
by the uniqueness of the weak limit for stochastic processes, it follows that $\mathcal Q (B^1(t)+x,B^2(t) +c\, t +x-\mathcal{E} )$ is equal in distribution to $B(t) + c \,t +x -\mathcal{E}$. \\

$2$) Define $Q^n$ as the length of the queue process constructed from the arrival $A^n$ and the service process $S$.
From Burke's theorem, we have that $\{ D^n(s): 0 \leq s < t\}$ is independent of $Q^n(t)$. On the other hand, the queue length process $Q$ constructed from the arrival $B^1(t)+x$ and the service $B^2(t) +c\, t +x-\mathcal{E}$ processes, 
has a explicit formulation in terms of the continuous operator $\mathcal{Q}$:
$$ Q(t) = B^1(t)+x - \mathcal{Q} (B^1(t)+x ,B^2(t) +c\, t +x-\mathcal{E} ).$$
So, by a similar argument to the one used to prove $1)$, we conclude that $Q^n$ converges weakly to $Q$ and the independence between $\{D^n(s) : 0 \leq s < t\}$ and $Q^n(t)$ is 
inherited by $D$ and $Q(t)$. \ep

\subsection{Convergence of the tandem Brownian queue}\label{Sec:MP}

We state a version of Borel-Cantelli lemma and prove a monotonicity result about the reflective operator that we will need. 

\begin{lema}\label{BC}
Let $(\Omega, \mathcal{F}, \pro)$ a probability space, $\{\mathcal{F}_n\}_{n \in \mathbb{N}}$ a filtration such that $\mathcal{F}= \bigcup_{n \in \na}\mathcal{F}_n$ and $O_n \in \mathcal{F}_n$, $n \in \na$. Then
$$ \{ O_n \, \, i.o. \} = \Big\{ \sum_{n \in \na} \pro (O_{n+1} | \mathcal{F}_{n} ) = \infty \Big\} \qquad a.s. .$$
\end{lema}
\bp See \cite{Kallenberg} pp. 108. \ep 

Let us denote by $\parallel \cdot \parallel _{[0,T]}$ the supremum norm on $[0,T]$.

\begin{lema}\label{SYNCH}
Denote the space of continuous real functions by $C[0,\infty)$. Let $f^1,f^2,g \in C[0,\infty)$ be such 
that $f^1(0),f^2(0) \geq g(0)$, and let $L_{f^1}(g)$, $L_{f^2}(g)$ be as in Equation \ref{birefl}. Then, for any $T>0$, 
$$  \parallel L_{f^1}(g) - L_{f^2}(g) \parallel_{[0,T]} \, \leq \,  \parallel f^1 - f^2 \parallel _{[0,T]} .$$
\end{lema}

\bp We have
\begin{eqnarray*}
 \parallel L_{f^1}(g) - L_{f^2}(g) \parallel_{[0,T]} & = &
\parallel \, [ g_t + \inf_{0 \leq s \leq t} \{ (f^1_s -g_s) \wedge 0 \} ] -  [ g_t + \inf_{0 \leq s \leq t} 
\{ (f^2_s - g_s) \wedge 0 \} ] \, \parallel _{[0,T]} \\
& = &  \parallel \sup_{0 \leq s \leq t} \{ (g_s-f^1_s) \vee 0 \} - \sup_{0 \leq s \leq t} \{ (g_s-f^2_s) \vee 0 \} \parallel _{[0,T]} \\
& \leq & \parallel (g_t-f^1_t) \vee 0  - (g_t-f^2_t) \vee 0 \parallel _{[0,T]} \\
& \leq & \parallel f^1 - f^2 \parallel _{[0,T]}.
\end{eqnarray*}
The first inequality follows from the Lipschitz continuity of the supremum mapping with Lipschitz constant equal to $1$ (see \cite{Whitt} pp. 436) and the second one holds since 
$ \parallel h^+ - i^+ \parallel_{[0,T]} \, \leq \, \parallel h-i \parallel_{[0,T]}$ for every pair $h,i$ of continuous real functions. \ep \\

\textbf{Proof of Theorem \ref{BROW:MP}:}

The proof relies on a coupling argument: we show that if different arrival processes are run through the same services, the resulting trajectories are eventually locally coupled. Since we know that there exists a stationary distribution under the reflecting dynamics of the tandem queues, given by the Theorem \ref{BROW:B}, we conclude the result. 
 
Fix $T>0$ and let $B^0$ be an independent Brownian motion. Define
$$ B^n = \Q \Big( B^{n-1}, W^n + c \, t - \sum_{i=1}^{n} \mathcal{E}^i  \Big),  \quad \forall n \geq 1 .$$
So, we apply the tandem queues dynamics to the process $B^0$ using the same service processes $ \{ W^n + c \, t \}_{n \in \mathbb{N}}$ and initial workloads $\{ \mathcal{E}^n \}_{n \in \mathbb{N}}$. By Theorem \ref{BROW:B}, $B^n + \sum_{i=1}^{n} \mathcal{E}^i $ has the law of a Brownian motion for every $n$, so in order to get convergence it is enough to prove that the trajectories of $A^n$ y $B^n$ eventually couple on $[0,T]$. 

The heart of the proof is the next: beginning with two different arrival process, at some step of the tandem dynamics we will have positive workload during a fixed period. Then the departures will coincide with the services in such period, and since we are using the same services processes, the departures of both systems will coincide. This coupling persists in the following iterations of the reflective dynamics. 

Last idea is resumed in the following observation. If 
$$W^{n+1}_t + c \, t -  \sum_{i=1}^{n+1} \mathcal{E}^i  \leq A^n_t, B^n_t, \quad \forall \, t \in [0,T],$$
then $A^{n+1}_t = B^{n+1}_t= W^{n}_t + c \, t -  \sum_{i=1}^{n+1} \mathcal{E}^i$, for $t \in [0,T]$, and the trajectories couple. 

For $n \in \mathbb{N}$, define the event
$$ O_n := \{ \, \omega \in \Omega: W^{n}_t + c \, t  - \sum_{i=1}^{n} \mathcal{E}^i \, \leq \, A_t ^{n-1}, B_t^{n-1} , \, \forall t \in [0,T] \, \},$$ 
and the $\sigma$-algebra $\mathcal{F}_{n} := \sigma( \{ A^0, B^0, W^i, \mathcal{E}^i : i \leq n \})$. Since $O_n$ is a coupling event
by Lemma \ref{BC} it is enough to prove that $\sum_{n=1}^{\infty} \esp ( 1_{O_{n+1}} | \mathcal{F}_n )= \infty$.

Define for $n \geq 0$, $\delta_n:=\parallel A^n_t -B^n_t \parallel_{[0,T]}$. Considering that the process $A_0$ does not explode in finite time a.s. we have that $\delta_0 < \infty$ a.s., hence:
$$ \sum_{n=1}^{\infty} \esp ( 1_{O_{n+1}} | \mathcal{F}_n ) = \sum_{k=1}^{ \infty} \sum_{n=1}^{\infty} 1_{ \{ k-1 \leq \delta_0 < k \} } \esp ( 1_{O_{n+1}} | \mathcal{F}_n ) = \sum_{k=1}^{ \infty} \sum_{n=1}^{\infty} \esp ( 1_{O_{n+1}}  1_{ \{ k-1 \leq \delta_0 < k \} } | \mathcal{F}_n ) ,$$
where $ 1_{ \{ k-1 \leq \delta_0 < k \} }$ is $\mathcal{F}_n$-measurable. We are using the same services and, by Lemma \ref{SYNCH}, we obtain $\delta_n( \omega) \leq \delta_0( \omega)$ for all $\omega \in \Omega$. Then
$$ \sum_{n=1}^{\infty} \esp ( 1_{O_{n+1}} | \mathcal{F}_n ) = \sum_{k=1}^{ \infty} \sum_{n=1}^{\infty} \esp ( 1_{O_{n+1}}   1_{ \{ k-1 \leq \delta_0 < k \} } 1_{ \{ \delta_n < k \} } | \mathcal{F}_n ). $$

Note now that $\{ W_t^{n+1} + ct - \sum_{i=1}^{n+1} \mathcal{E}^i \leq B^n - \delta_n, \, \forall t \in [0,T] \} \subseteq O_{n+1}$. Hence:
\begin{eqnarray*}
\sum_{n=1}^{\infty} \esp ( 1_{O_{n+1}} | \mathcal{F}_n ) & \geq & \sum_{k=1}^{ \infty} \sum_{n=1}^{\infty} \esp ( 1_{\{ W_t^{n+1} + ct- \sum_{i=1}^{n+1} \mathcal{E}^i \leq   \, B^n - \delta_n, \, \forall t \in [0,T] \}}   1_{ \{ k-1 \leq \delta_0 < k \} } 1_{ \{ \delta_n < k \} } | \mathcal{F}_n ) \\
& \geq & \sum_{k=1}^{ \infty} \sum_{n=1}^{\infty} \esp ( 1_{\{ W_t^{n+1} + ct - \sum_{i=1}^{n+1} \mathcal{E}^i \leq  \, B^n - k, \, \forall t \in [0,T] \}}   1_{ \{ k-1 \leq \delta_0 < k \} } 1_{ \{ \delta_n < k \} } | \mathcal{F}_n ) \\
& = & \sum_{k=1}^{ \infty} \sum_{n=1}^{\infty} \esp ( 1_{\{ W_t^{n+1} + ct - \sum_{i=1}^{n+1} \mathcal{E}^i \leq  \, B^n - k, \, \forall t \in [0,T] \}}   1_{ \{ k-1 \leq \delta_0 < k \} }  | \mathcal{F}_n ) \\
& = & \sum_{k=1}^{ \infty} \Big[ 1_{ \{ k-1 \leq \delta_0 < k \} } \sum_{n=1}^{\infty} \esp ( 1_{\{ W_t^{n+1} + ct - \sum_{i=1}^{n+1} \mathcal{E}^i \leq  \, B^n - k, \, \forall t \in [0,T] \}} | \mathcal{F}_n ) \Big].
\end{eqnarray*}

Define 
$$X_n^k:= \esp ( 1_{\{ W_t^{n+1} + ct -  \mathcal{E}^{n+1} < \, B^n + \sum_{i=1}^{n} \mathcal{E}^i - k, \, \forall t \in [0,T] \}} | \mathcal{F}_n ).$$

By Theorem \ref{BROW:B}, we have that $B^n + \sum_{i=1}^{n} \mathcal{E}^i$ is a Brownian motion for all $n$ and then the random variables $\{X_n^k \}_{n \in \mathbb{N}}$ are identically distributed. Moreover, since the reflective dynamics are Markovian with respect to the $n$-th step in the tandem queue, we have that 
$\{X_n^k \}_{n \in \mathbb{N}}$ are independent. By simple properties of the Brownian motion, they are strictly positive random variables. Therefore its sum $\sum_{n=1}^{\infty} X_n^k$ almost surely diverges for every $k$ and we are done by the 
almost sure finiteness of $\delta_0$. \ep





\bibliographystyle{abbrv}
\bibliography{BIBBROWN}

\end{document}